\documentclass{article}

\usepackage{algorithm, algpseudocode}
\usepackage{amsfonts, amsmath, amsopn, amsthm, epsfig}
  \numberwithin{equation}{section}
\usepackage{graphicx, epstopdf}
\usepackage{hyperref}
\usepackage{subfigure}

\usepackage{calc}
\usepackage[top=2.5cm,bottom=2.3cm,inner=2.3cm,outer=2.3cm,bindingoffset=0.5cm,footskip=0.55cm]{geometry}
\newtheorem{theorem}{Theorem}[section]
\newtheorem{corollary}[theorem]{Corollary}

\newtheorem{lemma}[theorem]{Lemma}

\theoremstyle{remark}
\newtheorem{remark}{Remark}

\newenvironment{lemma*}[2][Lemma]{\par\bgroup{\bfseries #1\ #2. }\it\ignorespaces}{\egroup}

\title{Online and Stochastic Universal Gradient Methods for Minimizing Regularized H\"older Continuous Finite Sums}

\author{Ziqiang Shi\footnotemark[1] \footnotemark[2], Rujie Liu\footnotemark[1]
}

\date{November 2013}

\newcommand{\BALD}{\begin{aligned}}
\newcommand{\EALD}{\end{aligned}}
\newcommand{\BALDS}{\begin{aligned*}}
\newcommand{\EALDS}{\end{aligned*}}
\newcommand{\BCAS}{\begin{cases}}
\newcommand{\ECAS}{\end{cases}}
\newcommand{\BEAS}{\begin{eqnarray*}}
\newcommand{\EEAS}{\end{eqnarray*}}
\newcommand{\BEQ}{\begin{equation}}
\newcommand{\EEQ}{\end{equation}}
\newcommand{\BIT}{\begin{itemize}}
\newcommand{\EIT}{\end{itemize}}
\newcommand{\BMAT}{\begin{bmatrix}}
\newcommand{\EMAT}{\end{bmatrix}}
\newcommand{\BNUM}{\begin{enumerate}}
\newcommand{\ENUM}{\end{enumerate}}

\newcommand{\BA}{\begin{array}}
\newcommand{\EA}{\end{array}}

\newcommand{\reals}{\mathbf{R}}

\DeclareMathOperator*{\minimize}{minimize}

\newcommand{\Prob}{\mathop{\mathbf{Pr}}}

\DeclareMathOperator{\sign}{sign}

\newcommand{\abs}[1]{\left| #1 \right|}

\newcommand{\norm}[1]{\left\| #1 \right\|}

\begin{document}

\maketitle

\renewcommand{\thefootnote}{\fnsymbol{footnote}}

\footnotetext[1]{Fujitsu Research \& Development Center, Beijing, China.}
\footnotetext[2]{shiziqiang@cn.fujitsu.com; shiziqiang7@gmail.com.}

\renewcommand{\thefootnote}{\arabic{footnote}}

\begin{abstract}
Online and stochastic gradient methods have emerged as potent tools in large scale optimization with both smooth convex and nonsmooth convex problems from the classes~$C^{1,1}(\reals^p)$ and~$C^{1,0}(\reals^p)$ respectively. However to our best knowledge, there is few paper to use incremental gradient methods to optimization the intermediate classes of convex problems with H\"older continuous functions~$C^{1,v}(\reals^p)$. In order fill the difference and gap between methods for smooth and nonsmooth problems, in this work, we propose the several online and stochastic universal gradient methods, that we do not need to know the actual degree of smoothness of the objective function in advance. We expanded the scope of the problems involved in machine learning to H\"older continuous functions and to propose a general family of first-order methods. Regret and convergent analysis shows that our methods enjoy strong theoretical guarantees. For the first time, we establish an algorithms that enjoys a linear convergence rate for convex functions that have H\"older continuous gradients.
\end{abstract}

\section{Introduction and problem statement}

Online and stochastic gradient methods (or referred to as incremental gradient methods) are of the most promising approaches in large scale machine learning tasks in these days~\cite{zinkevich2003online,wang2012online,mairal2013optimization,schmidt2013minimizing,xiao2014proximal,shalev2012proximal}.
Important advances of incremental gradient methods have been made on sequential learning in the recent literature on similar and famous problems, including lasso, logistic regression, ridge regression, and support vector regression. Composite objective mirror
descent (COMID)~\cite{duchi2010composite} generalizes mirror
descent~\cite{beck2003mirror} to the online setting.
Regularized dual averaging (RDA)~\cite{xiao2010dual} generalizes dual averaging~\cite{nesterov2009primal} to online and composite
optimization, and can be used for distributed optimization~\cite{duchi2012dual}. Online alternating direction multiplier method (ADMM)~\cite{suzuki2013dual}, RDA-ADMM~\cite{suzuki2013dual} and online proximal gradient (OPG) ADMM~\cite{wang2012online} generalize
classical ADMM~\cite{gabay1976dual} to online and stochastic settings. In stochastic gradient methods, more recent descent techniques like MISO~\cite{mairal2013optimization}, SAG~\cite{schmidt2013minimizing} and SVRG~\cite{xiao2014proximal} take update steps in the average gradient direction, and achieve linear convergence rate.

However, most current incremental gradient methods deal with smooth functions or non-smooth functions with Lipschitz-continues function values. In this paper, we consider incremental gradient methods with an objective function that has H\"older continuous gradients with degree $v$:
\begin{equation}
\|\nabla g(x) - \nabla g(y)\|_*  \leq M_v\|x - y\|^v,
  \label{eq:Hoelder-continuous-gradients}
\end{equation}
where $0 \leq v \leq 1$ and $\nabla g(x)$ means any subgradient if $g(x)$ is nonsmooth. It can be seen that $g(x)$ becomes smooth function with Lipschitz-continues gradients when $v=1$ and becomes non-smooth Lipschitz-continues function when $v=0$. $M_v$ is mainly used to characterize the variability of the (sub)gradients, all of this kind of functions form the class $C^{1,v}(\reals^p)$.
We consider the problems of the following form:
\begin{align}
  \minimize_{x \in \reals^p} \,f(x) := \frac{1}{n}\sum_{i=1}^n g_i (x) + h(x),
  \label{eq:composite-form-online}
\end{align}
where $g_i$ is a convex loss function with H\"older continuous gradients associated with a sample in a training set, and
$h$ is a convex penalty function or
regularizer. Let $g(x)=\frac{1}{n}\sum_{i=1}^n g_i (x)$.

If the Problem~(\ref{eq:composite-form-online}) is treated as minimizing of composite functions $g(x)+h(x)$, Nesterov has proposed the universal gradient methods (UGM) to solve it in~\cite{nesterov2013universal}. However, UGM for Problem~(\ref{eq:composite-form-online}) is a learning procedure in batch mode, which cannot deal with training data appearing in succession, such as audio processing~\cite{shi2013audio}. Furthermore, one can hardly ignore the fact that in reality the size of the data is rapidly increasing in
various domain and thus training set for the data probably cannot be loaded into the
memory simultaneously in batch mode methods.
In such situation, sequential learning becomes powerful tools.
In this paper, we generalize UGM to online and stochastic settings to deal with objective functions which have H\"older continuous gradients.

Assume $x^*$ is a solution of Problem (\ref{eq:composite-form-online}), and in this work, we introduce a novel kind of regret definition and seek bounds for this regret in the online learning setting with respect to $x^*$, defined as
\begin{equation}
        R(T,x^*,\epsilon):=\sum_{t=0}^T f_{g_t} (x_{t})- \sum_{t=0}^T f_{g_t}(x^*),
        \label{eq:online-regret}
\end{equation}
where $\epsilon$ if a pre-specified error limit. All of our algorithms need to first assume a fixed accuracy $\epsilon$, and then the smaller the $\epsilon$, the smaller the regret. For example, if we assume $\epsilon=1/T$, then we will have a regret bound of $O(1)$ after $T$ iterations. And if $\epsilon=1/\sqrt{T}$, then we will have a regret bound of $O(\sqrt{T})$ after $T$ iterations. Thus we have the results that look too good to be true, since our algorithms are different from previous online algorithms, and we have an extra parameter describing the accuracy. And the regret bound is not in a standard sense. Ours are in a sense that, for any fixed $T$, we can obtain an $O(1)$ bound after $T$ iterations.

We now outline the rest of the study. In Section~\ref{sec:ougm}, we propose online \textbf{prime}/\textbf{dual} universal gradient methods to solve the online optimization problem for the data that appear in succession and present the regret and convergence analysis. Section~\ref{sec:sug} states the stochastic universal gradient (SUG) method for the data that cannot be loaded into the
memory at the same time and show that the SUG achieves a linear convergence rate. We conclude in Section~\ref{sec:Conclusions}.

\subsection{Notations and lemmas}

Before proceeding, we introduce the notations and some useful lemmas formally first. In this work, we most adopt the nomenclature used by Nesterov on universal gradient methods~\cite{nesterov2013universal}. The functions encountered in this work are all convex if there are no other statements.

This inequality~\eqref{eq:Hoelder-continuous-gradients} ensures that
\begin{equation}
\label{eqn:holder-continuous-geometric}
   | g(x) -  g(y) - \nabla g(y)^T(x-y)| \leq \frac{M_v}{1+v} \|x-y\|^{1+v}.
\end{equation}

\emph{Bregman distance} is defined as
\begin{equation}
  \xi(x,y):=d(y)-d(x)-\langle\nabla d(x),y-x\rangle,
  \label{eq:bregman-divergence}
\end{equation}
where $d(x)$ is a \emph{prox-function}, which is differentiable strongly convex with convexity parameter equal to one and its minimum is $0$. Take derivative for $y$, we have
\[
    \nabla_y \xi(x,y) = \nabla d(y)-\nabla d(x).
\]

\emph{Bregman mapping} is defined as
\begin{equation}
  \hat{x} =\arg\min_{y}\bigl[g(x)+\langle\nabla g(x),y-x\rangle+M\xi(y,x)+h(y)\bigr],
  \label{eq:bregman-maping}
\end{equation}
where $h(y)$ is the fixed regularizer.

The first-order optimality condition for Problem~(\ref{eq:bregman-maping}) is
\begin{align}
  \langle\nabla g(x)+M(\nabla d(\hat{x})-\nabla d(x)) +\nabla h(\hat{x}), y-\hat{x}\rangle\geq 0.
  \label{eq:optimality-condition}
\end{align}

Some useful lemmas and equations introduced by~\cite{nesterov2013universal} are frequently employed in establishing the results and are stated below for the sake of completeness.
\begin{lemma}
  \label{lem:nesterov_lemma0}
If $\epsilon>0$ and $M>(\frac{1}{\epsilon})^{\frac{1-v}{1+v}}M_v^{\frac{2}{1+v}}$, then for any pair $t\geq 0$ we have
\begin{align}
  \frac{M_v}{1+v} t^{1+v}\leq \frac{1}{2}Mt^2+\frac{\epsilon}{2}.
  \label{eq:upgm-lemma0-1}
\end{align}
\end{lemma}

This lemma play an important role in this paper, which is been used to transform the  H\"older Continuous conditions to Lipschitz-continues conditions.

\begin{lemma}
  \label{lem:nesterov_lemma1}
If $g$ satisfy condition~\eqref{eq:Hoelder-continuous-gradients}, assume $\epsilon>0$ and $M>(\frac{1}{\epsilon})^{\frac{1-v}{1+v}}M_v^{\frac{2}{1+v}}$, then for any pair $x,y$ we have
\begin{align}
  g(y)\leq g(x)+\langle\nabla g(x),y-x\rangle+\frac{1}{2}M\|y-x\|^2+\frac{\epsilon}{2}.
  \label{eq:upgm-lemma1-1}
\end{align}
If $\hat{x}$ is the Bregman mapping at $x$ obtained by~\eqref{eq:bregman-maping}, then we have
\begin{align}
g(\hat{x})+h(\hat{x})\leq g(x)+\langle\nabla g(x),\hat{x}-x\rangle+M\xi(\hat{x},x)+h(\hat{x})+\frac{\epsilon}{2}.
  \label{eq:upgm-lemma1-2}
\end{align}
\end{lemma}

Throughout this work, we denote $\gamma(M_v,\epsilon):=(\frac{1}{\epsilon})^{\frac{1-v}{1+v}}M_\infty^{\frac{2}{1+v}}$.

\begin{lemma}
  \label{lem:nesterov_lemma2}
If $\phi(x)$ is convex and $\phi(x)-Md(x)$ is subdifferentiable, let $\bar{x}=\arg\min_{x}\phi(x)$, then we have
\begin{align}
      \phi(y)\geq \phi(\bar{x})+M\xi(\bar{x},y).
        \label{eq:lemma2-nesterov}
\end{align}
\end{lemma}

These lemmas are proposed in~\cite{nesterov2013universal},  please refer there for proofs if interested.

\section{Online Universal Gradient Method}
\label{sec:ougm}

In this section, we extend UGM to the online learning setting to deal with situation that the training data appearing in succession, such as multimedia information processing~\cite{shi2013audio}. The modification of UGM that we proposed is simple: just change $f_T(x)$ to $f_{g_t}(x)$ in each iteration and output the average value in each iteration. Our online algorithms are almost the same as the UGM with an important difference: we only meet and process one sample (one function) at each iteration. This methodology mainly comes from~\cite{duchi2010composite} and~\cite{wang2012online}. In the sequel, we consider three types of methods according to the original work of~\cite{nesterov2013universal}, from whose proofs we also draw some ingredients in ours.

\subsection{Online Universal Prime Gradient Method (O-UPGM)}

Lemma~\ref{lem:nesterov_lemma1} shows that the Bregman mapping can move the current point more close to the real solution, and this intuition form the core of the UGM and our online algorithms. In UGM, the Bregman mapping is employed to update the $x_t$ in each iteration, and $x_t$ is output as the solution after all the iterations.
Here we offer the general online universal primal gradient method (O-UPGM) solves Problem~(\ref{eq:composite-form-online}) in the following algorithm, where the same as UGM, Bregman mapping is also employed to update the $x_t$ in each iteration seeing current sample, while unlike UGM that the average of these $x_t$ is output as solutions after all the iterations.

\begin{algorithm}[H]
\caption{A generic O-UPGM}
\label{alg:General-online-UPGM}

\textbf{Input}: $L_0 > 0$ and $\epsilon > 0$.

1: \textbf{for} $t=0,1,\cdots,T$ \textbf{do}

2: Find the smallest $i_t \geq 0$ such that $g_t(\hat{x})+h(\hat{x})\leq g_t(x_t)+\langle\nabla g_t(x_t),\hat{x}-x_t\rangle+2^{i_t}L_t\xi(\hat{x},x_t)+h(\hat{x})+\frac{\epsilon}{2}$.

3: Set $x_{t+1}=\hat{x}$ and $L_{t+1}=2^{i_t-1}L_t$.

4: $t=t+1$.

5: \textbf{end for}

\textbf{Output}: $\bar{x}=\frac{1}{S_{T}}\sum_{t=1}^{T+1}\frac{1}{L_{t}}x_t$, where $S_{T} = \sum_{t=1}^{T+1}  \frac{1}{L_{t}}$.
\end{algorithm}

The above online UPGM is similar as UPGM except the $x_t$ update in O-UPGM uses a time varying function $f_{g_t}$. The following establishes the regret bound and the convergence rate for UPGM for general convex function with Hoelder continuous gradients.

\begin{theorem}
  \label{thm:general_o-upgm_regret_bound}
  Assume $M_v(g_t) < M_v$ and $h(x)$ is a simple convex function. Let the sequence $\{x_t\}$ be generated by the general O-UPGM in Algorithm~\ref{alg:General-online-UPGM}. Then we have
  \begin{align}
      \sum_{t=0}^{T} \frac{1}{L_{t+1}}[f_{g_t}(x_{t+1})-f_{g_t}(x^*)]\leq \frac{\epsilon}{2}S_{T}+2r_0(x^*),
       \label{eq:o-upgm-general-regret}
\end{align}
where $S_{T} = \sum_{t=1}^{T+1}  \frac{1}{L_{t}}$.
\end{theorem}

The ideas of the proof is near identical to that of UPGM by Nesterov~\cite{nesterov2013universal} but for completeness we give a simple version in the appendix.

We have the following remarks regarding the above result:
\begin{remark}
 All of our online algorithms (O-UPGM and O-UDGM) need to first assume a fixed accuracy $\epsilon$, and then the smaller the $\epsilon$, the more accurate the solution. For example, if we assume $\epsilon=1/T$, then we will have a regret bound of $O(1)$ after $T$ iterations. And if $\epsilon=1/\sqrt(T)$, then we will have a regret bound of $O(\sqrt(T))$ after T iterations. Thus we have the results that look too good to be true, since our algorithms are different from previous online algorithms, and we have an extra parameter describe the accuracy. And the regret bound is not in a standard sense. Ours are in a sense that, for any fixed $T$, we can obtain an $O(1)$ bound after T iteration.
\end{remark}

\begin{remark}
If we replace Step 2 and 3 in Algorithm~\ref{alg:General-online-UPGM} with $x_{t+1}=\mathfrak{B}_{2\gamma(M_v,\epsilon),g_t}(x_t)$, then $L_{t+1}=\gamma(M_v,\epsilon)$. Thus Theorem~\ref{thm:general_o-upgm_regret_bound} becomes
\begin{corollary}
\label{clry:upgm-regret-corollary}
  Assume $M_v(g_t) < M_v$ and $h(x)$ is a simple convex function. Let the sequence $\{x_t\}$ be generated by O-UPGM with fixed steps $L_{t+1}=\gamma(M_v,\epsilon)$. Then we have the standard regret bound
  \begin{align}
R(T,x^*,\epsilon) \leq \frac{\epsilon}{2}(T+1)+2r_0(x^*)\gamma(M_v,\epsilon).
  \end{align}
  Further, let $\epsilon=T^{-\frac{1+v}{2}}$, we have
\begin{align}
R(T,x^*,\epsilon)=O(T^{\frac{1-v}{2}}).
  \end{align}
\end{corollary}

\end{remark}

\subsection{Online Universal Dual Gradient Method (O-UDGM)}

The original UDGM is based on updating a simple model for objective function of Problem~(\ref{eq:composite-form-online}). We built a general online UDGM based on this principle for online or large scale problems.

\begin{algorithm}[H]
\caption{A generic O-UDGM}
\label{alg:General-online-UDGM}

\textbf{Input}: $L_0 > 0$, $\epsilon > 0$ and $\phi_0(x)=\xi(x_0,x)$.

1: \textbf{for} $t=0,1,\cdots,T$ \textbf{do}

2: Find the smallest $i_t \geq 0$ such that for point $x_{t,i_t}=\arg\min_{x} \phi_t(x)+\frac{1}{2^{i_t}L_t}[g_t(x_t)+\langle\nabla g_t(x_t),x-x_t\rangle+h(x)]$,  we have $f_{g_t}(\mathfrak{B}_{2^{i_t}L_t,g_t}(x_{t,i_t}))\leq \psi_{2^{i_t}L_t,g_t}^*(x_{t,i_t})+\frac{1}{2}\epsilon$.

3: Set $x_{t+1}=x_{t,i_t}$, $L_{t+1}=2^{i_t-1}L_t$ and $\phi_{t+1}(x)=\phi_t(x)+\frac{1}{2L_{t+1}}[g_t(x_t)+\langle\nabla g_t(x_t),x-x_t\rangle+h(x)]$.

4: $t=t+1$.

5: \textbf{end for}

\textbf{Output}: $\bar{x}=\frac{1}{S_{T}}\sum_{t=1}^{T+1}\frac{1}{L_{t}}x_t$, where $S_{T} = \sum_{t=1}^{T+1}  \frac{1}{L_{t}}$.
\end{algorithm}

\begin{theorem}
  \label{thm:general_o-udgm_regret_bound}
  Assume $M_v(g_t) < M_v$ and $h(x)$ is a simple convex function. Let the sequence $\{x_t\}$ be generated by the general O-UDGM. Then we have
 \begin{align}
      \sum_{t=0}^{T} \frac{1}{2L_{t+1}} f_{g_t}(x_t)- \sum_{t=0}^T \frac{1}{2L_{t+1}}f_{g_t}(x^*) \leq S_T\frac{\epsilon}{4}+\xi(x_0,x^*) \label{eq:o-udgm-general-regret}
\end{align}
  where $S_{T} = \sum_{t=1}^{T+1}  \frac{1}{L_{t}}$.
\end{theorem}

We have the following remarks regarding the above result:
\begin{remark}
If we replace Step 2 and 3 in Algorithm~\ref{alg:General-online-UDGM} with
 \begin{align}  \label{eq:dual-fixed-steps1}
 x_{t+1}=\arg\min_{x}\{\phi_t(x) +\frac{1}{2 \gamma(M_v,\epsilon)}[g_t(x_t)+\langle\nabla g_t(x_t),x-x_t\rangle+h(x)]\}
  \end{align}
and
 \begin{align}
\phi_{t+1}(x)=\phi_t(x)+\frac{1}{2 \gamma(M_v,\epsilon)}[g_t(x_t)+\langle\nabla g_t(x_t),x-x_t\rangle+h(x)]
\label{eq:dual-fixed-steps2}
 \end{align}
respectively,
then $L_{t+1}=\gamma(M_v,\epsilon)$ and Theorem~\ref{thm:general_o-udgm_regret_bound} becomes
\begin{corollary}
 \label{clry:udgm-regret-corollary}
  Assume $M_v(g_t) < M_v$ and $h(x)$ is a simple convex function. Let the sequence $\{x_t\}$ be generated by O-UDGM with fixed steps $L_{t+1}=\gamma(M_v,\epsilon)$. Then we have the standard regret bound
  \begin{align}
R(T,x^*,\epsilon) \leq \frac{\epsilon}{2}(T+1)+2\xi(x_0,x^*)\gamma(M_v,\epsilon).
  \end{align}
\end{corollary}

Further let $\epsilon=T^{-\frac{1+v}{2}}$, thus Corollary~\ref{clry:udgm-regret-corollary} becomes
\begin{corollary}
  Assume $M_v(g_t) < M_v$ and $h(x)$ is a simple convex function. Let the sequence $\{x_t\}$ be generated by the specific O-UDGM with $x_t$ updated by~(\ref{eq:dual-fixed-steps1}) and~(\ref{eq:dual-fixed-steps2}). Then we have
  \begin{align}
R(T,x^*,T^{-\frac{1+v}{2}})=O(T^{\frac{1-v}{2}}).
  \end{align}

\end{corollary}

\end{remark}

\section{Stochastic Universal Gradient Method}
\label{sec:sug}

In this section, we propose the stochastic universal gradient (SUG) method to deal with situation that the data probably cannot be loaded into the memory at the same time in batch mode methods since the size of the data is rapidly increasing.
We summarize the SUG method in Algorithm~\ref{alg:general-sug}.

\begin{algorithm}[H]
\caption{SUG: A generic stochastic universal gradient method}
\label{alg:general-sug}

\textbf{Input}: start point $x^0 \in$ dom $f$; for $i\in\{1,2,..,n\}$, let $g_i^0(x)=g_i(x^0)+(x-x^0)^T\nabla g_i(x^0)+M_0^i\xi(x^0,x)$, and $G^0(x)=\frac{1}{n}\sum_{i=1}^n g_i^0(x)$.

1: \textbf{repeat}

2: Solve the subproblem for new approximation of the solution: $x^{k+1} \leftarrow \arg\min_{x} \bigl[ G^k(x) + h(x)\bigr]$.

3: Sample $j$ from $\{1,2,..,n\}$, and update the surrogate functions:
\begin{align}
g_j^{k+1}(x)=g_j(x^{k+1})+(x-x^{k+1})^T\nabla g_j(x^{k+1})+M_{k+1}^i\xi(x^{k+1},x),
  \label{eq:subfunction-surrogate-update}
\end{align}
while leaving all other $g_i^{k+1}(x)$ unchanged: $g_i^{k+1}(x)\leftarrow g_i^{k}(x)$ ($i\neq j$); and $G^{k+1}(x)=\frac{1}{n}\sum_{i=1}^n g_i^{k+1}(x)$.

4: \textbf{until} stopping conditions are satisfied.

\textbf{Output}: $x^k$.
\end{algorithm}

\subsection{Convergence Analysis of SUG}

\begin{theorem}\label{thm:sug}
Suppose $g_i(x)$ satisfy condition~\eqref{eq:Hoelder-continuous-gradients} and $M \geq M_0^i> (\frac{2}{\epsilon})^{\frac{1-v}{1+v}}M_v^{\frac{2}{1+v}}$ for $i=1,...,n$, $d(x)$ satisfy $\|\nabla d(x) - \nabla d(y)\|_*  \leq M_d\|x - y\|^d$, $h(x)$ is strongly convex with $\mu_h\geq 0$, then the SUG iterations satisfy for $k\geq 1$:
\begin{equation}\label{eqn:converge-rate}
\mathbb{E}[f(x^k)] - f^* \le M \rho^{k-1} \|x^*-x^0\|^2+ \frac{3\epsilon}{4n\mu_h}\frac{1-\rho^{k-1}}{1-\rho}+\frac{3\epsilon}{4},
\end{equation}
where $\rho=\frac{1}{n} \frac{M}{\mu_h}+(1-\frac{1}{n})$.
\end{theorem}

We have the following remarks regarding the above result:
\begin{itemize}
\item
In order to satisfy $\mathbb{E}[f(x^k)] - f^* \leq \epsilon$,
the number of iterations~$k$ needs to satisfy
\[
    k\geq (\log\rho)^{-1} \log \bigl[\bigl(\frac{1}{4}-\frac{3}{4(\mu_h-M)}\bigr) \frac{\epsilon}{M\|x^*-x^0\|^2}\bigr]+1.
\]

\item
Inequality~\eqref{eqn:converge-rate} gives us a reliable stopping criterion for SUG method.

\end{itemize}

Since $\mathbb{E}[f(x^k)] - f^* \geq 0$, Markov's
inequality and Theorem~\ref{thm:sug} imply
that for any $\epsilon>0$,
\[
\Prob\Bigl( f(x^k) - f^* \geq \epsilon\Bigr)
~\leq~ \frac{\mathbb{E}[f(x^k)] - f^*}{\epsilon}
~\leq~ \frac{M \rho^{k-1} \|x^*-x^0\|^2}{\epsilon}+ \frac{3}{4n\mu_h}\frac{1}{1-\rho}+\frac{3}{4}.
\]
Thus we have the following high-probability bound.

\begin{corollary}\label{cor:sug-high-prob-1}
Suppose the assumptions in Theorem~\ref{thm:sug} hold.
Then for any $\epsilon>0$ and $\delta\in(0,1)$, we have
\[
  \Prob\bigl(f(x^k) - f(x^\star) \leq \epsilon \bigr) \geq 1-\delta
\]
provided that the number of iterations~$k$ satisfies
\[
k \geq (\log\rho)^{-1} \log \bigl[\bigl(\delta-\frac{3}{4}-\frac{3}{4(\mu_h-M)}\bigr) \frac{\epsilon}{M\|x^*-x^0\|^2}\bigr]+1.
\]
\end{corollary}

\section{Conclusions}
\label{sec:Conclusions}

In this paper, in order to fill the difference and gap between methods for smooth and nonsmooth problems, we propose efficient online and stochastic gradient algorithms to optimization the intermediate classes of convex problems with H\"older continuous functions~$C^{1,v}(\reals^p)$. We establish regret bounds for the objective and linear convergence rates for convex functions that have Hoelder continuous gradients. There are some directions that the current study can be extended. In this paper, we have focused on the theory; it would be meaningful to also do the numerical evaluation and implementation details, and we give some simple applications in Section~\ref{sec:appendix_application}. Second, combine with randomized block coordinate method~\cite{nesterov2012efficiency} for minimizing regularized convex functions with a huge number of varialbes/coordinates. Moreover, due to the trends and needs of big data, we are designing distributed/parallel SUG for real life applications. In a broader context, we believe that the current paper could serve as a basis for examining the method for the classes of convex problems with H\"older continuous functions~$C^{1,v}(\reals^p)$.
%Section~\ref{sec:appendix_application}

\appendix
\section*{Appendix}

In this Appendix, we give the proofs of the propositions.

\section{Proof of Theorem~\ref{thm:general_o-upgm_regret_bound}}

First, we show that the algorithm, especially step 2 is well defined. Due to~(\ref{eq:upgm-lemma1-1}) and~(\ref{eq:upgm-lemma1-2}) in Lemma~\ref{lem:nesterov_lemma1} and monotonically increasing of $2^{i_t}L_t$, when $L_{t+1}>\gamma(M_v,\epsilon)$, we have $f_{g_t}(\mathfrak{B}_{2^{i_t}L_t,g_t}(x_t))\leq \psi_{2^{i_t}L_t,g_t}^*(x_t)+\frac{1}{2}\epsilon$.
Thus we always have
\[
        2L_{t+1}=2^{i_t}L_t \leq 2\gamma(M_v,\epsilon).
\]
Let us fix an arbitrary point $y$, and denote $r_t(y):=\xi(x_t,y)$. Then we have
\begin{eqnarray*}
        r_{t+1}(y)&=&d(y)-d(x_{t+1})-\langle\nabla d(x_{t+1}),y-x_{t+1}\rangle   \\
        &\leq&  d(y)-d(x_{t+1})-\langle\nabla d(x_t),y-x_{t+1}\rangle+\frac{1}{2L_{t+1}}\langle\nabla g_t(x_t)+\nabla h(x_{t+1}),y-x_{t+1}\rangle
\end{eqnarray*}
and
\begin{eqnarray*}
      &  &d(y)-d(x_{t+1})-\langle\nabla d(x_t),y-x_{t+1}\rangle  \\
       &\leq  &d(y)-d(x_t)-\langle\nabla d(x_t),x_{t+1}-x_{t}\rangle-\frac{1}{2}\|x_{t+1}-x_{t}\|^2-\langle\nabla d(x_t),y-x_{t+1}\rangle \\
        &= & r_t(y)-\frac{1}{2}\|x_{t+1}-x_{t}\|^2,
\end{eqnarray*}
where the first inequality is derived using~(\ref{eq:optimality-condition}) and the second inequality is obtained using the strong convex property of the prox-function $d(x)$.

Thus we have
\begin{eqnarray*}
       & &r_{t+1}(y)-r_t(y)  \\  &\leq &\frac{1}{2L_{t+1}}\langle\nabla g_t(x_t)+\nabla h(x_{t+1}),y-x_{t+1}\rangle-\frac{1}{2}\|x_{t+1}-x_{t}\|^2  \\
        &= &\frac{1}{2L_{t+1}}\langle\nabla h(x_{t+1}),y-x_{t+1}\rangle +\frac{1}{2L_{t+1}}\langle\nabla g_t(x_t),y-x_t\rangle  -\frac{1}{2L_{t+1}}(\langle\nabla g_t(x_t),x_{t+1}-x_t\rangle+L_{t+1}\|x_{t+1}-x_{t}\|^2)   \\
         &\leq &\frac{1}{2L_{t+1}}[h(y)-h(x_{t+1})+g_t(x_t)-g_t(x_{t+1})+\frac{\epsilon}{2}+\langle\nabla g_t(x_t),y-x_t \rangle]
\end{eqnarray*}

Thus we obtain
\begin{eqnarray*}
       \frac{1}{2L_{t+1}}f_{g_t}({x_{t+1}})+r_{t+1}(y)  \leq \frac{1}{2L_{t+1}}[g_t(x_t)+\langle\nabla g_t(x_t),y-x_t \rangle+h(y)+\frac{\epsilon}{2}]+r_t(y).
\end{eqnarray*}

Summing up, we have
\[
      \sum_{t=0}^{T} \frac{1}{L_{t+1}}f_{g_t}({x_{t+1}})+r_{T+1}(y)  \leq \sum_{t=0}^{T}  \frac{1}{L_{t+1}}[g_t(x_t)+\langle\nabla g_t(x_t),y-x_t \rangle+h(y)+\frac{\epsilon}{2}]+2r_0(y).
\]

Let $y=x^*$, we have $g_t(x_t)+\langle\nabla g_t(x_t),x^*-x_t \rangle\leq g_t(x^*)$, thus
\[
      \sum_{t=0}^{T} \frac{1}{L_{t+1}}[f_{g_t}(x_{t+1})-f_{g_t}(x^*)]\leq \frac{\epsilon}{2}\sum_{t=0}^{T}  \frac{1}{L_{t+1}}+2r_0(x^*)
\]
and it is proved.

\section{Proof of Theorem~\ref{thm:general_o-udgm_regret_bound}}

Similar with the reasoning of Theorem~\ref{thm:general_o-upgm_regret_bound}, that the algorithm, especially step 2 is well defined, and we also always have
\[
        2L_{t+1}=2^{i_t}L_t \leq 2\gamma(M_v,\epsilon).
\]

Denote $y_t=\mathfrak{B}_{2^{i_t}L_t,g_t}(x_t)$ and $\phi_t^*=\arg\min_{x}\phi_t(x)$. Let $S_t=\sum_{i=0}^{t} \frac{1}{L_{i+1}}$, we first prove that
\begin{align}
       \sum_{i=0}^{t} \frac{1}{2L_{i+1}} f_{g_i}(y_i) \leq \phi_{t+1}^* + S_t\frac{\epsilon}{4}
        \label{eq:dual-target-bound}
\end{align}
is valid for all $t\geq 0$. Indeed, for $t=0$ we have
\[
        f_{g_0}(y_0)-\frac{\epsilon}{2} \leq \psi_{2^{i_0} L_0,g_0}^*(x_0)
       = g_0(x_0)+\langle\nabla g_0(x_0),y_0-x_0\rangle+2^{i_0} L_0\xi(x_0,y_0)+h(y_0)   =  2^{i_0} L_0\phi_{1}^*.
\]

In view of~(\ref{eq:lemma2-nesterov}) in Lemma~\ref{lem:nesterov_lemma2}, for any $t\geq 0$, we have
\[
       \phi_{t+1}(x)\geq \phi_{t+1}(x_t)+\xi(x_t,x).
\]

Assume that~(\ref{eq:dual-target-bound}) is true for some $t\geq 0$. Then
\begin{eqnarray*}
       & &\min_x\phi_{t+2}(x)  \\ &\geq &\min_x \{\phi_{t+1}(x)+\frac{1}{2 L_{t+2}}[g_{t+1}(x_{t+1})   +\langle\nabla g_{t+1}(x_{t+1}),x-x_{t+1}\rangle+h(x)]\}  \\
        &\geq &\min_x \{\phi_{t+1}(x_{t+1})+\xi(x_{t+1},x)+\frac{1}{2 L_{t+2}}[g_{t+1}(x_{t+1})+\langle\nabla g_{t+1}(x_{t+1}),x-x_{t+1}\rangle+h(x)]\} \\
       &\geq &  \phi_{t+1}(x_{t+1})+ \frac{1}{2 L_{t+2}}[f_{g_{t+1}}(y_{t+1})-\frac{\epsilon}{2}]  \geq  -S_{T}\frac{\epsilon}{4}+\sum_{i=0}^{t+1} \frac{1}{2L_{i+1}} f_{g_i}(y_i).
\end{eqnarray*}
Thus~(\ref{eq:dual-target-bound}) is proved.

In view of~(\ref{eq:dual-target-bound}), we have
\[
       \sum_{i=0}^{t} f_{g_i}(y_i) \leq 2\gamma(M_v,\epsilon)\phi_{t+1}^* + \frac{\epsilon}{2}(t+1).
\]
Since
\begin{eqnarray*}
    \phi_{t+1}(y)&\leq &\phi_{t}(y) + \frac{1}{2L_{t+1}}[g_t(y)+h(y)] \leq \phi_{t-1}(y) + \frac{1}{2L_{t}}f_{g_{t-1}}(y) +\frac{1}{2L_{t+1}}f_{g_t}(y)  \\
    &\leq & \sum_{i=0}^t \frac{1}{2L_{i+1}}f_{g_i}(y)+\xi(x_0,y),
\end{eqnarray*}
we have
\[
       \sum_{t=0}^T \frac{1}{2L_{t+1}}f_{g_t}(y)+\xi(x_0,y) \geq -S_T\frac{\epsilon}{4}+\sum_{t=0}^{T} \frac{1}{2L_{t+1}} f_{g_t}(y_t).
\]
Rearrange the terms, and let $y=x^*$ the theorem is proved.

\section{Proof of Theorem~\ref{thm:sug}}

Since in each iteration of the SUG, we obtain a function $g_i^{k}(x)$ with random parameters to approximate each $g_i(x)$:
\begin{align}
g_i^{k}(x)=g_i(x^{\theta_{i,k}})+(x-x^{\theta_{i,k}})^T\nabla g_i(x^{\theta_{i,k}})+M_{\theta_{i,k}}^i\xi(x^{\theta_{i,k}},x),
  \label{eq:subfunction-surrogate-assum}
\end{align}
where $\theta_{i,k}$ is a random variable which have the following conditional probability distribution in each iteration:
\[
\mathbb{P}(\theta_{i,k}=k|j)=\frac{1}{n} \quad \text{and} \quad \mathbb{P}(\theta_{i,k}=\theta_{i,k-1}|j)=1-\frac{1}{n},
\]
that yields
\begin{align}
\label{eq:x-update-relation}
\mathbb{E}[\|x^*-x^{\theta_{i,k}}\|^2]=\frac{1}{n}\mathbb{E}[\|x^*-x^k\|^2]+(1-\frac{1}{n})\mathbb{E}[\|x^*-x^{\theta_{i,k-1}}\|^2].
\end{align}

Since $M_{\theta_{i,k}}^i > (\frac{2}{\epsilon})^{\frac{1-v}{1+v}}M_v^{\frac{2}{1+v}}$, by lemma~\ref{lem:nesterov_lemma1} we have
\[
    g_i(x)\leq g_i(x^{\theta_{i,k}})+(x-x^{\theta_{i,k}})^T\nabla g_i(x^{\theta_{i,k}})+M_{\theta_{i,k}}^i\xi(x^{\theta_{i,k}},x)+\frac{\epsilon}{4}.
\]
Thus by~\eqref{eq:subfunction-surrogate-assum}, we have
\[
    g_i(x)\leq g_i^k(x)+\frac{\epsilon}{4}
\]
and summing over $i=1,\ldots,n$ yields
\begin{align}
\label{eq:nesterov-bound}
    g(x)\leq G^k(x)+\frac{\epsilon}{4}.
\end{align}
Take derivative of~\eqref{eq:subfunction-surrogate-assum}, we have
\[
   \|\nabla g_i^k(x) - \nabla g_i^k(y)\| = \|\nabla d(x) - \nabla d(y)\|.
\]

Set $\delta_i^k(x)=g_i(x)-g_i^k(x)$, we have
\begin{eqnarray*}
   |\delta_i^k(x) -\delta_i^k(y) -\langle \nabla \delta_i^k(y),x-y\rangle | &=& \| \int_0^1 \langle\nabla\delta_i^k(y+t(x-y))  -\nabla\delta_i^k(y),x-y\rangle  {\rm d}t\| \\
   &\le& \int_0^1\| \langle\nabla\delta_i^k(y+t(x-y))  -\nabla\delta_i^k(y),x-y\rangle \| {\rm d}t \\
    &\le&  \int_0^1 \|\nabla\delta_i^k(y+t(x-y))  -\nabla\delta_i^k(y)\|\|x-y \| {\rm d}t \\
    &\le&  \|y-x \|\int_0^1 \|\nabla g_i^k(y+t(x-y)) - \nabla g_i^k(y)\| {\rm d}t \\
    &+&  \|y-x\|\int_0^1 \|\nabla g_i(y+t(x-y)) - \nabla g_i(y)\| {\rm d}t \\
    &\le&  \|y-x \|\int_0^1 t^dM_d\|x - y\|^d {\rm d}t + \|y-x\|\int_0^1 t^vM_v\|x - y\|^v {\rm d}t\\
    &\le&  \frac{1}{1+d}M_d\|y-x \|^{1+d} + \frac{1}{1+v}M_v\|y-x \|^{1+v}\\
\end{eqnarray*}

Let $y=x^{\theta_{i,k}}$, and since we have $\delta_i^k(x^{\theta_{i,k}})=0$ and $\nabla\delta_i^k(x^{\theta_{i,k}})=0$, thus
\[
   |g_i(x)-g_i^k(x)| \le  \frac{1}{1+d}M_d\|x-x^{\theta_{i,k}} \|^{1+d} + \frac{1}{1+v}M_v\|x-x^{\theta_{i,k}} \|^{1+v}.
\]

Summing over $i=1,\ldots,n$ yields
\begin{align}
\label{eq:diff_bound}
    [G^{k}(x)+h(x)] - [g(x)+h(x)] \leq \frac{1}{n}\sum_{i=1}^n  \bigl[ \frac{1}{1+d}M_d\|x-x^{\theta_{i,k}} \|^{1+d} + \frac{1}{1+v}M_v\|x-x^{\theta_{i,k}} \|^{1+v} \bigr].
\end{align}

Since $G^{k}(x)+h(x)$ is $\mu_h$-strongly convex, together with~\eqref{eq:diff_bound} and~\eqref{eq:nesterov-bound}, we have
\begin{eqnarray*}
f(x^{k+1})+\frac{\mu_h}{2}\|x-x^{k+1}\|^2
&\leq& G^{k}(x^{k+1})+h(x^{k+1})+\frac{\epsilon}{4}+\frac{\mu_h}{2}\|x-x^{k+1}\|^2 \\
&\leq& G^{k}(x)+h(x)+\frac{\epsilon}{4} \\
&=& f(x)+[G^{k}(x)+h(x)-f(x)]+\frac{\epsilon}{4} \\
&\leq& f(x)+ \frac{1}{n}\sum_{i=1}^n  \bigl[ \frac{1}{1+d}M_d\|x-x^{\theta_{i,k}} \|^{1+d} + \frac{1}{1+v}M_v\|x-x^{\theta_{i,k}} \|^{1+v} \bigr]+\frac{\epsilon}{4} .
\end{eqnarray*}
By taking the expectation of both sides and let $x=x^*$, using~\eqref{eq:upgm-lemma0-1} yields
\begin{eqnarray*}
   \mathbb{E}[f(x^{k+1})] - f^* &\leq& \mathbb{E}[\frac{1}{n}\sum_{i=1}^n  \bigl[ \frac{1}{1+d}M_d\|x^*-x^{\theta_{i,k}} \|^{1+d} + \frac{1}{1+v}M_v\|x^*-x^{\theta_{i,k}} \|^{1+v} \bigr]]-\mathbb{E}[\frac{\mu_h}{2}\|x^*-x^{k+1}\|^2]+\frac{\epsilon}{4} \\
   &\leq& \mathbb{E}[\frac{M}{n}\sum_{i=1}^n  \bigl[ \|x^*-x^{\theta_{i,k}} \|^2 \bigr]]-\mathbb{E}[\frac{\mu_h}{2}\|x^*-x^{k+1}\|^2]+\frac{3\epsilon}{4}.
\end{eqnarray*}

Thus we have
\[
    \mu_h \|x^{k+1}-x^*\|^2  \leq  \mathbb{E}[\frac{M}{n}\sum_{i=1}^n  \bigl[ \|x^*-x^{\theta_{i,k}} \|^2 \bigr]]+\frac{3\epsilon}{4}.
\]
then we have
\begin{eqnarray*}
  \mathbb{E}[\frac{1}{n}\sum_{i=1}^n \|x^*-x^{\theta_{i,k}}\|^2]
  &=&  \frac{1}{n}\|x^k-x^*\|^2+(1-\frac{1}{n})\mathbb{E}[\frac{1}{n}\sum_{i=1}^n \|x^*-x^{\theta_{i,k-1}}\|^2] \\
   &\leq&  [\frac{1}{n} \frac{M}{\mu_h}+(1-\frac{1}{n})]\mathbb{E}[\frac{1}{n}\sum_{i=1}^n \|x^*-x^{\theta_{i,k-1}}\|^2]+ \frac{3\epsilon}{4n\mu_h}\\
    &\leq&  [\frac{1}{n} \frac{M}{\mu_h}+(1-\frac{1}{n})]^k\mathbb{E}[\frac{1}{n}\sum_{i=1}^n \|x^*-x^{\theta_{i,0}}\|^2] \\&+& \frac{3\epsilon}{4n\mu_h}\bigl(1+(\frac{1}{n} \frac{M}{\mu_h}+(1-\frac{1}{n}))+...+(\frac{1}{n} \frac{M}{\mu_h}+(1-\frac{1}{n})(^{k-1}\bigr)\\
      &\leq&  [\frac{1}{n} \frac{M}{\mu_h}+(1-\frac{1}{n})]^k \|x^*-x^0\|^2+ \frac{3\epsilon}{4n\mu_h}\frac{1-(\frac{1}{n} \frac{M}{\mu_h}+(1-\frac{1}{n}))^{k}}{1-\frac{1}{n} \frac{M}{\mu_h}+(1-\frac{1}{n})}.
\end{eqnarray*}

Thus we have
\[
\mathbb{E}[f(x^{k+1})] - f^* \le M [\frac{1}{n} \frac{M}{\mu_h}+(1-\frac{1}{n})]^k \|x^*-x^0\|^2+ \frac{3\epsilon}{4n\mu_h}\frac{1-(\frac{1}{n} \frac{M}{\mu_h}+(1-\frac{1}{n}))^{k}}{1-\frac{1}{n} \frac{M}{\mu_h}+(1-\frac{1}{n})}+\frac{3\epsilon}{4}.
\]

\section{Some applications}
\label{sec:appendix_application}

In this appendix, we present some applications of our methods.

\subsection{Online lasso problem}
The lasso problem is formulated as follows:
\begin{align}
  \minimize_{x \in \reals^{p}}\,\frac{1}{n}\sum_{t=1}^n \|a_t^Tx-b_t\|^2 + \mu\norm{x}_1,
  \label{eq:lasso}
\end{align}
where $a_t, x\in \reals^{p}$ and $b_t$ is a scalar.

Throughout this section, let $g(x)=\frac{1}{n}\sum_{t=1}^n \|a_t^Tx-b_t\|^2$ and $h(x)=\mu\norm{x}_1$, $d(x)=\frac{1}{2}\|x\|^2$, then $\xi(x,y)=\frac{1}{2}\|x-y\|^2$. The subdifferential of the function $\|a^Tx-b\|^2$ is $2(a^Tx-b)a$.

The Bregman mapping associate with $g(x)$ and the component function $g_t(x)=\|a_t^Tx-b_t\|^2$ are
\begin{align}
  \hat{x}&=\arg\min_{y} \{\frac{1}{T}\sum_{t=1}^T \|a_t^Tx-b_t\|^2+\langle\frac{2}{T}\sum_{i=1}^T (a_t^Tx-b_t)a_t,y-x\rangle+M\frac{1}{2}\|x-y\|^2  + \mu\norm{y}_1\} \nonumber \\ &
  =\sign(x-\frac{2}{MT}\sum_{i=1}^T (a_t^Tx-b_t)a_t)\cdot\max\{\abs{x-\frac{2}{TM}\sum_{i=1}^T (a_t^Tx-b_t)a_t}-\frac{\mu}{M},0\}
  \label{eq:lasso-bregman-maping}
\end{align}
and
\begin{align}
  \hat{x}&=\arg\min_{y} \{ \|a_t^Tx-b_t\|^2+\langle\ 2(a_t^Tx-b_t)a_t,y-x\rangle+M\frac{1}{2}\|x-y\|^2  + \mu\norm{y}_1\} \nonumber \\ &
  =\sign(x-\frac{2}{M} (a_t^Tx-b_t)a_t)\cdot\max\{\abs{x-\frac{2}{M} (a_t^Tx-b_t)a_t}-\frac{\mu}{M},0\}
  \label{eq:online-lasso-bregman-maping}
\end{align}
respectively.

In online UDGM and SUG, we have
\begin{eqnarray*}
\phi_{t+1}(x)&=\phi_{t}(x)+a_t[g_t(x_t)+\langle\nabla g_{t}(x_t),x-x_t\rangle+\mu\norm{x}_1]
\\&= \xi(x_0,x) + \sum_{i=1}^{t} a_i[g_{i}(x_{i})+\langle\nabla g_{i}(x_i),x-x_i\rangle+\mu\norm{x}_1].
\end{eqnarray*}
Then we have
\begin{eqnarray*}
x_{t+1}&=\arg\min_{x}\phi_{t+1}(x) = \arg\min_{x} \{\frac{1}{2}\|x_0-x\|^2 +  \sum_{i=1}^{t} a_i[\langle\nabla g_{i}(x_i),x\rangle+\mu\norm{x}_1] \}  \\ &= \sign(x_0-\sum_{i=1}^{t} a_i\nabla g_{i}(x_i))\cdot\max\{\abs{x_0-\sum_{i=1}^{t} a_i\nabla g_{i}(x_i)}-\mu\sum_{i=1}^{t} a_i,0\}.
\end{eqnarray*}

\subsection{Online Steiner problem}

In continuous Steiner problem we are given by centers $c_i\in \reals^{p}$, $i=1,...,m$. It is necessary to find the optimal location of the service center $x$, which minimizes the total distance to all other centers. Thus, our problem is as follows:
\begin{align}
  \min_{x \in \reals^{p}} \,g(x) := \frac{1}{m}\sum_{i=1}^m \|x-c_i\|,
  \label{eq:online-steiner}
\end{align}
where all norms in this problem are Euclidean. UGM solves that problem effectively. However, in real application, new locations will be added to the system, such as new shop opening or new warehouse establishing. Thus our online and stochastic gradient algorithms are needed.

Let $h(x)=0$, $d(x)=\frac{1}{2}\|x\|^2$, then $\xi(x,y)=\frac{1}{2}\|x-y\|^2$. The subdifferential of the Euclidean norm $\|x\|$ is $\frac{x}{\|x\|}$ if $x\neq 0$ or $\{g|\|x\|\leq 1\}$ if $x=0$. In order to simplify the formula, we here denote
$\nabla \|x\|=  \frac{x}{\|x\|}$ instead distinguishing between $x=0$ and $x\neq 0$.

The Bregman mapping associate with $\frac{1}{m}\sum_{i=1}^m \|x-c_i\|$ and the component function $\|x-c_i\|$ are
\begin{align}
  \hat{x}=&\arg\min_{y} \{\frac{1}{m}\sum_{i=1}^m \|x-c_i\|+\langle\frac{1}{m}\sum_{i=1}^m \frac{x-c_i}{\|x-c_i\|},y-x\rangle+M\frac{1}{2}\|x-y\|^2  \}
  =x-\frac{1}{mM}\sum_{i=1}^m\frac{x-c_i}{\|x-c_i\|}
  \label{eq:steiner-bregman-maping}
\end{align}
and
\begin{align}
  \hat{x}=\arg\min_{y} \{ \|x-c_i\|+\langle \frac{x-c_i}{\|x-c_i\|},y-x\rangle+M\frac{1}{2}\|x-y\|^2  \}
  =x-\frac{1}{M}\frac{x-c_i}{\|x-c_i\|}
  \label{eq:online-steiner-bregman-maping}
\end{align}
respectively

In online UDGM and SUG for Steiner problem, we have
\begin{eqnarray*}
\phi_{t+1}(x)=\phi_t(x)+a_t[g_{t}(x_t)+\langle\nabla g_{t}(x_t),x-x_t\rangle]
= \xi(x_0,x) + \sum_{i=1}^{t} a_i[g_{i}(x_{i})+\langle\nabla g_{i}(x_i),x-x_i\rangle]
\end{eqnarray*}
where $ g_{i}(x_{i})= \|x_{i}-c_{i}\|$ and $\nabla g_{i}(x_i)= \frac{x_{i}-c_{i}}{\|x_{i}-c_{i}\|}$.
Thus we have
\begin{eqnarray*}
x_{t+1}&=\arg\min_{x}\phi_{t+1}(x) = \arg\min_{x} \frac{1}{2}\|x_0-x\|^2 + \sum_{i=1}^{t} a_{i}\langle \frac{x_{i}-c_{i}}{\|x_{i}-c_{i}\|},x\rangle \\ &= \arg\min_{x} \frac{1}{2}\|x_0-x\|^2 + \langle \sum_{i=1}^{t}a_{i}\frac{x_{i}-c_{i}}{\|x_{i}-c_{i}\|},x\rangle =  x_0 - \sum_{i=1}^{t}a_{i}\frac{x_{i}-c_{i}}{\|x_{i}-c_{i}\|}.
\end{eqnarray*}

\bibliography{OnlineUGM}
\bibliographystyle{plain}

\end{document}